\documentclass[12pt]{article}
\usepackage{amsmath}
\usepackage{amssymb}
\usepackage{dsfont}
\usepackage{cite}
\newsymbol \blackbox 1004
\newcommand{\eh}{\hfill}\newlength{\sperr}

\newenvironment{proof}{{\settowidth{\sperr}{\bf\rm
Proof}%
\par\addvspace{0.3cm}\noindent\parbox[t]{1.3\sperr}
{\bf\rm P\eh r\eh o\eh o\eh f\eh }%
}}{\nopagebreak\mbox{}
$\blackbox$\par\addvspace{0.3cm}}

\def\nn{\nonumber}

\def\a{\alpha}

\def\Lam{\Lambda}
\def\s{\sigma}
\def\la{\lambda}
\def\om{\omega}

\def\wh{\widehat}
\def\wt{\widetilde}
\def\wc{\check}
\def\ov{\overline}

\def\p{\partial}

\def\BC{{\mathbb C}}

\def\BR{{\mathbb R}}
\def\BN{{\mathbb N}}

\def\cla{{\mathcal A}}

\def\cld{{\mathcal D}}

\def\clh{{\mathcal H}}

\def\clj{{\mathcal J}}

\def\clu{{\mathcal U}}

\def\clr{\mathcal{R}}
\def\cls{\mathcal{S}}

\def\clu{{\mathcal U}}

\def\cld{{\mathcal D}}

\def\tr{{\rm tr}}
\newcommand{\E}{\mathrm{e}}
\newcommand{\I}{\mathrm{i}}

\newtheorem{Pa}{Paper}[section]
\newtheorem{Tm}[Pa]{{\bf Theorem}}

\newtheorem{Cy}[Pa]{{\bf Corollary}}
\newtheorem{Rk}[Pa]{{\bf Remark}}

\newtheorem{Ee}[Pa]{{\bf Example}}

\newtheorem{Pn}[Pa]{{\bf Proposition}}

\title{Einstein, $\s$-model and Ernst-type equations \\
and non-isospectral GBDT version \\ of  Darboux transformation}

\author{Alexander Sakhnovich}

\date{}
\begin{document}
\maketitle

\begin{abstract} We  present a non-isospectral GBDT version of B\"acklund-Darboux transformation 
for the gravitational and $\sigma$-model equations. New families of explicit solutions correspond  
to the case of GBDT with non-diagonal generalized matrix eigenvalues.
An interesting integrable Ernst-type system, the auxiliary linear systems of which are
non-isospectral canonical systems, is studied as well.
\end{abstract}

{MSC(2010): 35Q51, 35Q76, 83C15, 37J05

\vspace{0.2em}

{\bf Keywords:} Einstein equation, $\sigma$-model, Ernst-type equation, Darboux transformation,
generalized matrix eigenvalue, non-isospectral canonical system.

\section{Introduction} \label{intro}
\setcounter{equation}{0}
The study of the integrable reductions of Einstein field equations goes back to the seminal
paper \cite{BeZ} (see also \cite{Mai}). The  survey \cite{Alev2} includes  several references to the interesting
articles which precede \cite{BeZ} and a bibliography of the related works during  thirty years after its publication.
For the recent references
one can turn, for instance, to \cite{LeM}. Following the publication of \cite{BeZ}, a closely related $\s$-model equation was studied
in \cite{MiYa}.  Gravitational (Einstein)  equation and $\s$-model equation both belong to the so called
non-isospectral case where the spectral parameter depends on other variables (see, e.g., \cite{BZM, Calo} on this topic).
We apply to the gravitational (Einstein) equation and $\s$-model equation the non-isospectral GBDT version of B\"acklund-Darboux
transformation. This version of B\"acklund-Darboux transformation is especially suitable for the explicit construction of the wave functions
and solutions of those equations. Generalized
matrix eigenvalues $\cla$ are used in GBDT instead of the  usual eigenvalues, and {\it new classes of explicit solutions} appear when we deal with the non-diagonal 
$\cla$ (e.g., $\cla$ in the normal Jordan form).

Hamiltonian evolution equations are related to Einstein and $\s$-model equations (and play an essential role in their study),
see, for instance, \cite{Ash0, Ash, Iri}. In this paper, we investigate an interesting Ernst-type integrable nonlinear system:
\begin{align}\nn &
H(\xi,\eta)-\clh(\xi,\eta)=\I \big( \clh(\xi,\eta)J H(\xi,\eta)-H(\xi,\eta)J \clh(\xi,\eta)\big), \\
\nn &
H_{\eta}(\xi,\eta)=\clh_{\xi}(\xi,\eta) \quad (H\geq 0, \quad \clh \geq 0), \quad J=\begin{bmatrix}0 & I_p \\ I_p & 0\end{bmatrix},
\end{align}
where the Hamiltonians $H$ and $\clh$ are $2p \times 2p$ matrix functions and the auxiliary (to our Ernst-type system) linear systems are non-isospectral canonical
systems with these Hamiltonians. As usual, $\I$ above stands for the  imaginary unit ($\I^2=-1$) and $I_p$ is the $p\times p$ identity matrix.

Gravitational (Einstein) equation in light-cone coordinates has the form
\begin{align} &       \label{1.1}
\big(\a(\xi,\eta)u_{\xi}(\xi,\eta) u(\xi,\eta)^{-1})_{\eta}+\big(\a(\xi,\eta)u_{\eta}(\xi,\eta) u(\xi,\eta)^{-1})_{\xi}=0, \quad \a_{\xi\eta}=0,
 \end{align} 
where $\a$ is a scalar function,  $u$ is a $2\times 2$ matrix function,
 and $u_{\xi}=\frac{\p}{\p \xi} u$. Physically meaningful solutions $u$ of \eqref{1.1} have the properties \cite{BeZ, BZM}:
\begin{align} &       \label{1.2}
\a\in \BR, \quad u(\xi,\eta)\in {\mathrm{GL}}(2, \BR),
 \end{align} 
where $\BR$ is the real axis and ${\mathrm{GL}}(2, \BR)$ stands for the set of $2\times 2$ invertible matrix functions
with real-valued entries.
The solutions satisfying an additional property 
\begin{align} &       \label{3.4+}
\det(u)=\a^2.
 \end{align} 
are constructed via the multiplication of the solutions $u$ of \eqref{1.1} satisfying \eqref{1.2}
by certain real-valued scalar functions (see \eqref{3.4} or \cite[(2.17)]{BeZ}).

In the important paper \cite{MiYa}, the authors wrote down $\s$-model equation in the form \eqref{1.1},
where $u$ are $m\times m$ invertible matrix functions with complex-valued entries ($u\in {\mathrm{GL}}(m, \BC))$.
More precisely, it is supposed that the  relations
\begin{align} &       \label{1.3}
\a(\xi,\eta)\in \BR, \quad u(\xi,\eta)\in {\mathrm{GL}}(m, \BC), \quad u(\xi,\eta)^*Ju(\xi,\eta) \equiv J
 \end{align} 
hold.  Here, $\BC$ stands for the complex plane, $u^*$ means complex conjugate transpose of $u$, and we assume further that the $m\times m$ matrix $J$
satisfies relations
\begin{align} &       \label{1.4}
J=J^*=J^{-1}.
 \end{align} 
 
The paper consists of five sections. Some basic GBDT relations for the equation \eqref{1.1} are given in Section \ref{Prel}. Using these relations, we 
express (in Section \ref{grav}) wide families of solutions of the   $\s$-model and gravitational equations (so called {\it transformed} solutions)
via some initial solutions. The discussed above Ernst-type equation is studied in Section \ref{Ernst}.
Finally, explicit formulas and new explicit solutions are presented in Section \ref{Some}.
Several auxiliary results are proved in the Appendices \ref{apA} and \ref{apB}.

Some notations have been introduced above and further notations are explained here.
The notation $\BN$ means the set of positive integer numbers, 
$\ov{\a}$ stands for the complex conjugate of $\a$, and  the inequality $H\geq 0$ for some matrix $H$ means that $H=H^*$ 
and has nonnegative eigenvalues. 
The set of $i \times k$ matrices with real-valued
entries in denoted by $\BR^{i\times k}$. The spectrum of matrix $\cla$ is denoted by $\s(\cla)$.
We say that the function is continuously differentiable if its first derivatives exist and are
continuous (in the topology $\BR^k$ if it is a function of $k$ variables).
\section{Preliminaries} \label{Prel}
\setcounter{equation}{0}
{\bf 1.} 
B\"acklund-Darboux transformations and related commutation methods
present an important tool in spectral,
gauge  and soliton theories (see, e.g., \cite{Ci, D, Ge, GeT, Gu, KoSaTe, Mar, MS, Mi, ZM}).
Our GBDT version of B\"acklund--Darboux transformation was first introduced in \cite{SaA1}
(see further results and references in the papers \cite{KaaSa, KoSaTe, SaA2, ALS01, ALS20} and in the
book \cite{SaSaR}).

In this section, we derive some important relations for the  equation
\eqref{1.1} from our more general GBDT results in \cite[Sections 2,3]{SaA2}. 
We study the case \eqref{1.1}, \eqref{1.3}. We discuss also the modification
of the solution of \eqref{1.1}, \eqref{1.2} such that \eqref{3.4+} holds.

Integrable linear equations are often considered in the so called zero curvature form 
\begin{align} &       \label{2.1}
\frac{d}{d \eta}G-\frac{d}{d \xi}F+[G,F]=0 \quad ([G,F]:=GF-FG),
 \end{align} 
 which is the compatibility condition
 of the auxiliary linear systems
\begin{equation}       \label{2.2}
\frac{d}{d \xi}w(\xi,\eta,\la)=G(\xi,\eta,\la)w(\xi,\eta,\la), \quad \frac{d}{d \eta}w(\xi,\eta,\la)=F(\xi,\eta,\la)w(\xi,\eta,\la),
 \end{equation} 
where $\xi$ and $\eta$ are independent variables,
 $\la$ is the  spectral parameter and $w$ is an $m\times m$ non-degenerate matrix function
(fundamental solution). In the so called  isospectral case, where $\la$ does not depend on $\xi$ and $\eta$,
one can write, for instance, $G_{\eta}=\frac{\p}{\p \eta}G$ instead of $\frac{d}{d \eta}G$. In the non-isospectral case,
where $\la$ depends on $\xi$ and (or) $\eta$, we need the total derivatives $\frac{d}{d \xi}$ and (or) $\frac{d}{d \eta}$
with respect to these variables.

Here, relation \eqref{2.1} easily follows from \eqref{2.2}
but the fact that \eqref{2.1} yields the existence  of $w$ satisfying \eqref{2.2} is somewhat more complicated, see \cite{ALS12}
and references therein. Clearly, zero curvature representation \cite{AKNS, ZM0, TF}  is closely related to Lax pairs.

According to \cite{BZM}, equation \eqref{1.1} is equivalent to 
\eqref{2.1} in the case
 \begin{align} &       \label{2.3}
 G(\xi,\eta,\la)=-\frac{1}{\la -1}q(\xi,\eta), \quad F(\xi,\eta,\la)=-\frac{1}{\la +1}Q(\xi,\eta).
 \end{align} 
 Moreover, the case is non-isospectral, that is, $\la$ is a scalar function depending on
 the variables
 $\xi$ and $\eta$ and on the ``hidden spectral parameter" $z$. The dependence
 of $\la$ on $\xi$ and $\eta$ is given by the equations \cite{BZM, SaA2}:
 \begin{align} &       \label{2.4}
\la_{\xi}=-\frac{\a_{\xi}}{\a}\la \frac{\la+1}{\la-1}=-\frac{\a_{\xi}}{\a}\la-\frac{2\a_{\xi}}{\a}- \frac{2\a_{\xi}}{\a(\la-1)},
\\ &       \label{2.5}
\la_{\eta}=-\frac{\a_{\eta}}{\a}\la \frac{\la-1}{\la+1}=-\frac{\a_{\eta}}{\a}\la+\frac{2\a_{\eta}}{\a}- \frac{2\a_{\eta}}{\a(\la+1)}.
 \end{align}
 The equality $a_{\xi\eta}=0$ in \eqref{1.1} means that $\a$ admits representation:
 \begin{align} &       \label{2.6}
\a(\xi, \eta)=f(\xi)+h(\eta).
 \end{align}
Since $\a=\ov{\a}$ (see \eqref{1.2} and \eqref{1.3}), we assume further that
\begin{align} &       \label{2.6'}
f(\xi)=\ov{f(\xi)}, \quad h(\eta)=\ov{h(\eta)}
 \end{align}
in \eqref{2.6}. 
It follows from \eqref{2.4}--\eqref{2.6} (see \cite{BZM}) that one can choose
  \begin{align} &       \label{3.1}
\la(\xi, \eta,z)=\frac{h(\eta)-f(\xi)-z+\sqrt{(z-2h(\eta))(z+2f(\xi))}}{f(\xi)+h(\eta)}.
 \end{align}
 \begin{Rk}\label{RkR} The functions $f, h$ and the branch of the square root in \eqref{3.1} $($or, equivalently, in \eqref{3.8}$)$
 should be chosen so that $\la(\xi,\eta)$ is well-defined  and
continuously differentiable. For this purpose, we may also either restrict the domains of $\xi$ and $\eta$
or turn to the Riemann surfaces (see, e.g., \cite[p. 510]{MiYa}).
\end{Rk}
We note  that the matrix functions $q$ and  $Q$ in \eqref{2.3} are connected with the solution $u$ of the corresponding
equation  \eqref{1.1} by the equalities 
 \begin{align} &       \label{2.7}
q(\xi,\eta)=u_{\xi}(\xi,\eta) u(\xi,\eta)^{-1}, \quad Q(\xi,\eta)=-u_{\eta}(\xi,\eta) u(\xi,\eta)^{-1}
\end{align}
(see \cite{BZM} or  \cite[(44)]{SaA2}).
\begin{Rk}\label{RkComp}
 In fact, taking into account the property  \eqref{3.9} of $\la$, we rewrite \eqref{2.1}  in the form
\begin{align} &\label{C}
q_{\eta}+Q_{\xi}=[q,Q], \quad (\a q)_{\eta}=(\a Q)_{\xi}.
\end{align}
The existence of $u$ satisfying \eqref{2.7} and the fact that  \eqref{1.1} holds for this $u$ easily follow from \eqref{C}. \end{Rk}

{\bf 2.}  In this paper, each generalized B\"acklund-Darboux transformation (GBDT) is determined  by 
some  initial system \eqref{2.2}, \eqref{2.3} (to which GBDT is applied)
and by a triple of matrices $\{\cla, S(0,0), \Pi(0,0)\}$,
 where $\cla$ and $S(0,0)$ are $n\times n$  matrices $(n\in \BN)$, $\Pi(0,0)$ is an $n\times m$ matrix and the matrix identity
 \begin{align} &       \label{2.8}
\cla S(0,0)-S(0,0)\cla^*=\I\Pi(0,0)J\Pi(0,0)^*, \quad J=J^*=J^{-1}.
 \end{align}
 holds.  Clearly, instead of the initial system  \eqref{2.2}, \eqref{2.3}, we may fix 
 the functions $\a(\xi,\eta), \, q(\xi, \eta)$ and $Q(\xi,\eta)$ generating \eqref{2.2}, \eqref{2.3}
 or the functions
 $\a(\xi,\eta)$ and $u(\xi,\eta)$ satisfying \eqref{1.1} (in which case $q$ and $Q$
are given by  \eqref{2.7}).
 
 We assume that 
 \begin{align} &       \label{2.14}
\a=\ov{\a}, \quad S(0,0)=S(0,0)^*, \quad qJ=-J q^*, \quad QJ=-J Q^*.
 \end{align}
 
 Below, we show that similar to the isospectral case, the so called Darboux matrix function $w_A$  has at each $\xi$
 and $\eta$ the form of the  transfer matrix function:
 \begin{align} &       \label{2.8+}
w_A(\xi,\eta,\la)=I_m-\I J\Pi(\xi,\eta)^*S(\xi, \eta)^{-1}\big(A(\xi,\eta)-\la I_n\big)^{-1}\Pi(\xi, \eta).
 \end{align}
The transfer matrix function was introduced in this form by Lev Sakhnovich in  \cite{SaL1} (see also \cite{SaL2, SaSaR}).
The corresponding matrix functions $A(\xi,\eta), \, \Pi(\xi,\eta)$ and $S(\xi,\eta)$ in \eqref{2.8+} are defined by the values
 $A(0,0)=\cla$, $\Pi(0,0)$ and $S(0,0)$, respectively, and by the linear equations:
\begin{align} &       \label{2.9}
A_{\xi}=-\frac{\a_{\xi}}{\a}A-\frac{2\a_{\xi}}{\a}I_n- \frac{2\a_{\xi}}{\a}(A-I_n)^{-1},
\\ &       \label{2.10}
A_{\eta}=-\frac{\a_{\eta}}{\a}A+\frac{2\a_{\eta}}{\a}I_n- \frac{2\a_{\eta}}{\a}(A+I_n)^{-1},
\\ &       \label{2.11}
\Pi_{\xi}=(A-I_n)^{-1}\Pi q, \quad \Pi_{\eta}=(A+I_n)^{-1}\Pi Q,
\\ &       \label{2.12}
S_{\xi}=\frac{\a_{\xi}}{\a}\left(S-2(A-I_n)^{-1}S(A^*-I_n)^{-1}\right)-\I(A-I_n)^{-1}\Pi qJ\Pi^*(A^*-I_n)^{-1},
\\ &       \label{2.13}
S_{\eta}=\frac{\a_{\eta}}{\a}\left(S-2(A+I_n)^{-1}S(A^*+I_n)^{-1}\right)-\I(A+I_n)^{-1}\Pi QJ\Pi^*(A^*+I_n)^{-1}.
 \end{align}
 Recall that  (for our non-isospectral case) $\la =\la(\xi,\eta, z)$ in \eqref{2.8+}. The Darboux matrix function $w_A$ transforms the fundamental
 solution $w$ of the initial system \eqref{2.2} into the fundamental solution $w_A w$ of the transformed system.
 
In view of \eqref{2.9} and \eqref{2.10} we have 
$$A_{\xi \eta}=A_{\eta \xi}=2\frac{\a_{\xi}\a_{\eta}}{\a^2}A^3(A-I_n)^{-1}(A+I_n)^{-1},$$
and so the compatibility condition for systems \eqref{2.9}, \eqref{2.10} is fulfilled.
In order to see that equations \eqref{2.11} are compatible, we take into account \eqref{2.9}, \eqref{2.10} and differentiate $\Pi_{\xi}$ with respect
 to $\eta$ and $\Pi_{\eta}$ with respect to $\xi$. It follows that
 \begin{align}& \nn
 \a \Pi_{\xi \eta}=(A-I_n)^{-1}(A+I_n)^{-1}\big(\a_{\eta}A\Pi q+\a \Pi Qq+\a(A+I_n)\Pi q_{\eta}\big), \\
 & \nn
 \a \Pi_{\eta \xi}=(A-I_n)^{-1}(A+I_n)^{-1}\big(\a_{\xi}A\Pi Q+\a \Pi qQ+\a(A-I_n)\Pi Q_{\xi}\big).
 \end{align}
Now, the compatibility condition $ \Pi_{\xi \eta}= \Pi_{\eta \xi}$ is immediate from \eqref{C}.
The equality $S_{\xi \eta}=S_{\eta \xi}$ is proved in a similar way although more complicated
calculations are required for that purpose. See some further details in Appendix~\ref{apA}.

Equations \eqref{2.9}--\eqref{2.11} are derived from
the more general formulas considered in \cite[pp. 1252-1254]{SaA2}. The equality $\a=\ov{\a}$ enabled us to set (for our special case) $A_1=A$ and $A_2=A^*$
 in the formula \cite[(19)]{SaA2}.  After substitution $\Pi_2(0,0)^*=\I J\Pi_(0,0)^*$, formula \cite[(6)]{SaA2} at the point $(0,0)$ took the form
 \eqref{2.8}. Next, we  used the last two equalities in \eqref{2.14} in order to set $\Pi_1\equiv \Pi$ and $\Pi_2^*\equiv \I J\Pi^*$ in \cite[(5),\,(19),\,(20))]{SaA2}.
 Equations \eqref{2.9}--\eqref{2.11} followed. Now, from \eqref{2.8}--\eqref{2.13} one obtains (see \cite[(6)]{SaA2}):
\begin{align} &       \label{2.15}
A(\xi,\eta)S(\xi,\eta)-S(\xi,\eta)A(\xi,\eta)^*=\I \Pi(\xi,\eta)J\Pi(\xi,\eta)^*.
 \end{align} 
 Finally, formulas (11)--(13) in \cite{SaA2} imply that
  \begin{align}   \label{2.15+}&    \frac{d}{d \xi}w_A(\xi, \eta, \la)=\wh G(\xi, \eta, \la)w_A(\xi, \eta, \la)
-w_A(\xi, \eta, \la)G(\xi, \eta, \la),
\\ \label{2.16}&
\frac{d}{d \eta}w_A(\xi, \eta, \la)=\wh F(\xi, \eta, \la)w_A(\xi, \eta, \la)
-w_A(\xi, \eta, \la)F(\xi, \eta, \la),
\\  \label{2.17}&
\wh G(\xi,\eta,\la)=-\frac{1}{\la -1}\wh q(\xi,\eta), \quad F(\xi,\eta,\la)=-\frac{1}{\la +1}\wh Q(\xi,\eta).
 \end{align}
Here, according to \cite[(10),\,(13)]{SaA2} we have the following expressions for the transformed coefficients $\wh q$ and $\wh Q$
(denoted by $\wh q_{11}$ and $\wh Q_{11}$ in \cite{SaA2}):
\begin{align}        \nn
\wh q=&\big(I_m-\I J\Pi^*S^{-1}(A-I_n)^{-1}\Pi \big)q\big(I_m+\I J\Pi^*(A^*-I_n)^{-1}S^{-1}\Pi \big)
\\ \label{2.18} &
-2\I (\a_{\xi}/\a)J\Pi^*S^{-1}(A-I_n)^{-1}S(A^*-I_n)^{-1}S^{-1}\Pi,
\\   \nn
\wh Q=&\big(I_m-\I J\Pi^*S^{-1}(A+I_n)^{-1}\Pi \big)Q\big(I_m+\I J\Pi^*(A^*+I_n)^{-1}S^{-1}\Pi \big)
\\ \label{2.19} &
-2\I (\a_{\eta}/\a)J\Pi^*S^{-1}(A+I_n)^{-1}S(A^*+I_n)^{-1}S^{-1}\Pi.
 \end{align}
 When we invert $S$ above, we consider the corresponding 
 formulas in the points of invertibility of $S$. Note that equalities \eqref{2.12}--\eqref{2.14} yield
 $S(\xi,\eta)=S(\xi,\eta)^*$, and so $ \wh q J$ and $ \wh Q J$ given by \eqref{2.18} and \eqref{2.19} satisfy skew-self-adjointness conditions similar
 to the last two equalities in \eqref{2.14} for $  q J$ and $  Q J$:
 \begin{align} &       \label{2.20}
\wh qJ=-J({\wh q})^*, \quad \wh QJ=-J (\wh Q)^*.
 \end{align}
 Recall that $\a$ and $u$ satisfy \eqref{1.1}, and so \eqref{2.1} holds. Hence, according to  \cite[Theorem 6.1]{SaSaR}
the initial system \eqref{2.2} is compatible. 
 \begin{Rk}\label{Comp1} For simplicity, we assume in the text that   $G$ and $F$ are continuously differentiable 
 and, for this purpose, $\a(\xi,\eta)=f(\xi)+h(\eta),$ $q(\xi,\eta)$ and $Q(\xi,\eta)$ are  continuously differentiable
$($or, instead of the requirements on $q$ and $Q$, that $u$ is two times continuously differentiable$)$.
In fact,  the conditions which we need in order that \eqref{2.1} yields the compatibility
of systems \eqref{2.2} are weaker $($see, e.g., \cite[Theorem 6.1]{SaSaR}$)$.
\end{Rk}
In view of \eqref{2.2} and taking into account  \eqref{2.15+}, \eqref{2.16} 
we see that the matrix function $\wh w(\xi, \eta, \la)=w_A(\xi, \eta, \la)w(\xi, \eta, \la)$ satisfies the system
 \begin{align} &       \label{2.2'}
\frac{d}{d \xi}\wh w(\xi,\eta,\la)=\wh G(\xi,\eta,\la)\wh w(\xi,\eta,\la), \quad \frac{d}{d \eta}\wh  w(\xi,\eta,\la)=\wh F(\xi,\eta,\la)\wh w(\xi,\eta,\la).
 \end{align} 
Thus, the transformed system \eqref{2.2'} is compatible and the compatibility condition holds:
\begin{align} &       \label{2.1'}
\frac{d}{d \eta}\wh G-\frac{d}{d \xi} \wh F+[\wh G,\wh F]=0 .
 \end{align}
 Moreover, relations \eqref{2.2'} (or, equivalently, relations \eqref{2.15+} and \eqref{2.16}) show that
 $w_A$ of the form \eqref{2.8+} is, indeed, a Darboux   matrix function.
 
{\bf 3.}  Finally, we describe a way to modify a solution  of \eqref{1.1} (when $m=2$) so that
the modified solution $u$ satisfies the equality \eqref{3.4+}
It is easy to see that  a $2\times 2$ matrix function $u=\{u_{ik}\}_{i,k=1}^2$ has a property:
 \begin{align} &       \label{3.1+}
u_{\xi}u^{-1}=\frac{1}{\det (u)}\begin{bmatrix} (u_{11})_{\xi}u_{22}-(u_{12})_{\xi}u_{21} & * \\ * &  (u_{22})_{\xi}u_{11}-(u_{21})_{\xi}u_{12}
\end{bmatrix}.
 \end{align}
 Thus, considering traces ``$\tr$" of both sides of \eqref{3.1+} we have
 \begin{align} &       \label{3.2}
\tr\big(u_{\xi}u^{-1}\big)=\frac{(\det u)_{\xi}}{\det(u)}.
 \end{align}
Clearly, a similar to \eqref{3.2} formula is valid for $\tr\big(u_{\eta}u^{-1}\big)$. Hence, taking traces in \eqref{1.1}  one obtains
\begin{align} &       \label{3.3}
\left(\a \frac{(\det u)_{\xi}}{\det(u)}\right)_{\eta}+\left(\a \frac{(\det u)_{\eta}}{\det(u)}\right)_{\xi}=0.
 \end{align}
Now, assuming that  $\a$ and   some $2\times 2$ matrix function $\wc u$ satisfy \eqref{1.1}, \eqref{1.2}, and that $\det (\wc u)>0$, 
one (using standard calculations) derives that $\a$ and the matrix function 
 \begin{align} &       \label{3.4}
u:=\a (\det \wc u)^{-1/2}\wc u
 \end{align}
 satisfy \eqref{1.1}, \eqref{1.2} and equality \eqref{3.4+}.

 \section{$\sigma$-model and gravitational equations} \label{grav}
\setcounter{equation}{0}
{\bf 1.}  Clearly, \eqref{1.1} remains valid if we multiply $u$ (from the right) by some constant $m\times m$ matrix.
Hence, without loss of generality one may assume that 
 \begin{align} &       \label{3.4!}
u(0,0)=I_m.
 \end{align}
First, we prove   the following theorem on the construction of solutions of the $\s$-model equation.
\begin{Tm} \label{TmSM} Let $\a$ and $u$ satisfy   equation \eqref{1.1},
let a triple of matrices $\{\cla, S(0,0), \Pi(0,0)\}$ satisfying \eqref{2.8}
be given and assume that relations  \eqref{2.14} and \eqref{3.4!} hold,
where $q$ and $Q$ in \eqref{2.14} are given by \eqref{2.7}. Set
\begin{align} &       \label{3.7}
\clu(\xi,\eta):=I_m-\I J\Pi(\xi,\eta)^*S(\xi,\eta)^{-1}A(\xi,\eta)^{-1}\Pi(\xi,\eta),
 \end{align}
where  the matrix functions $A(\xi,\eta), \, \Pi(\xi,\eta)$ and $S(\xi,\eta)$ are introduced by the linear equations  \eqref{2.9}--\eqref{2.13}.

Then, the scalar function $\a$ and the matrix function
\begin{align} &       \label{3.7+}
\wh u(\xi,\eta)=\clu(\xi,\eta)u(\xi, \eta) 
 \end{align}
satisfy equation \eqref{1.1}, that is,
\begin{align} &       \label{1.1'}
\big(\a(\xi,\eta)\wh u_{\xi}(\xi,\eta) \wh u(\xi,\eta)^{-1})_{\eta}+\big(\a(\xi,\eta)\wh u_{\eta}(\xi,\eta)\wh u(\xi,\eta)^{-1})_{\xi}=0.
 \end{align} 
Moreover, $\wh u$ is $J$-unitary:
\begin{align} &       \label{1.2'}
\wh u(\xi,\eta)^*J\wh u(\xi,\eta)\equiv J.
\end{align} 
\end{Tm}
\begin{proof}. Fixing the branch of the square root we rewrite \eqref{3.1} in the form
  \begin{align} &       \label{3.8}
\la(\xi, \eta,z)=\frac{\big(\sqrt{z-2h(\eta)}-\sqrt{z+2f(\xi)}\big)^2}{z-2h(\eta)-(z+2f(\xi))}=\frac{\sqrt{z-2h(\eta)}-\sqrt{z+2f(\xi)}}{\sqrt{z-2h(\eta)}+\sqrt{z+2f(\xi)}}.
 \end{align}
Thus, we have
\begin{align} &       \label{3.9}
\la \to 0 \quad {\mathrm{for}} \quad z\to \infty.
 \end{align}
 Relations \eqref{2.4}, \eqref{2.5} and \eqref{3.9} yield
 \begin{align} &       \label{3.10}
\la_{\xi} \to 0, \quad \la_{\eta}\to 0 \quad {\mathrm{for}} \quad z\to \infty.
 \end{align}
For a fixed constant $\mu=\la(\xi,\eta,z_{\mu})$, from the definition \eqref{2.8+} we obtain
 \begin{align} \label{3.11} &    \big(I_m-\I J\Pi(\xi,\eta)^*S(\xi,\eta)^{-1}(A(\xi,\eta)-\mu I_n)^{-1}\Pi(\xi,\eta)\big)_{\xi}
\\  \nn & 
 =\big(w_A(\xi, \eta, \la)\big)_{\xi}
 +\I \la_{\xi}(\xi,\eta,z_{\mu})J\Pi(\xi,\eta)^*S(\xi,\eta)^{-1}(A(\xi,\eta)-\mu I_n)^{-2}\Pi(\xi,\eta);
\\  \label{3.12} &    \big(I_m-\I J\Pi(\xi,\eta)^*S(\xi,\eta)^{-1}(A(\xi,\eta)-\mu I_n)^{-1}\Pi(\xi,\eta)\big)_{\eta}
\\  \nn & 
 =\big(w_A(\xi, \eta, \la)\big)_{\eta}
 +\I\la_{\eta}(\xi,\eta,z_{\mu})J\Pi(\xi,\eta)^*S(\xi,\eta)^{-1}(A(\xi,\eta)-\mu I_n)^{-2}\Pi(\xi,\eta).
 \end{align}
 When $z_{\mu}$ tends to infinity, formulas \eqref{2.15+}--\eqref{2.17} and \eqref{3.9}--\eqref{3.12}
 yield
 \begin{align} &       \label{3.13}
\clu_{\xi}(\xi,\eta)=\wh q(\xi,\eta)\clu(\xi,\eta)- \clu(\xi,\eta)q(\xi,\eta), \\
&       \label{3.14}
\clu_{\eta}(\xi,\eta)=-\wh Q(\xi,\eta)\clu(\xi,\eta)+ \clu(\xi,\eta)Q(\xi,\eta).
 \end{align} 
Using equalities \eqref{2.7},  \eqref{3.13} and \eqref{3.14} we derive
 \begin{align} &       \label{3.15}
\wh u_{\xi}(\xi,\eta)= \wh q(\xi,\eta) \wh u(\xi,\eta), \quad \wh u_{\eta}(\xi,\eta)=-\wh Q(\xi,\eta) \wh u(\xi,\eta),
\end{align} 
where $\wh q$, $\wh Q$ and $\wh u$ are given by \eqref{2.18}, \eqref{2.19} and \eqref{3.7+}, respectively.
Applying Remark \ref{RkComp} to the zero curvature equation \eqref{2.1'}  and taking into account \eqref{3.15}, we see
that $\wh u$ satisfies  equation \eqref{1.1'}.
Moreover,  relations \eqref{1.4}, \eqref{2.20} and \eqref{3.15} imply that
\begin{align} &       \label{3.16}
\big(\wh u(\xi,\eta)^*J\wh u(\xi,\eta)\big)_{\xi}=0,  \quad \big(\wh u(\xi,\eta)^*J\wh u(\xi,\eta)\big)_{\eta}=0.
 \end{align}
From \eqref{2.8+} and \eqref{2.15} one obtains 
 \begin{align} &       \label{E19}
w_A(\xi,\eta,\ov{\la})^*J\
w_A(\xi,\eta,\la) \equiv J
 \end{align}  
(see, e.g., \cite{SaL1} or \cite[(1.84)]{SaSaR}). In particular, we have
\begin{align} &       \label{3.17}
\clu(\xi,\eta)^*J\clu(\xi,\eta)=J.
 \end{align}
 By virtue of \eqref{3.4!}, \eqref{3.7+} and \eqref{3.17}, the equality 
 \begin{align} &       \label{3.18}
\wh u(0,0)^*J\wh u(0,0)=J
 \end{align}
holds.  Finally, formulas \eqref{3.16} and \eqref{3.18} yield \eqref{1.2'}.
\end{proof}

{\bf 2.}  Setting in Theorem \ref{TmSM} $m=2$, we easily obtain the following corollary for the gravitational equation.
\begin{Cy} \label{Cygr} Let $\a$ and $u$ satisfy \eqref{1.1} and \eqref{1.2}, let a triple of matrices $\{\cla, S(0,0), \Pi(0,0)\}$,
which satisfies the matrix identity \eqref{2.8}, be given
and assume that the relations  
 \begin{align} &       \label{3.19}
\cla, \, S(0,0)\in {\mathrm{GL}}(n, \BR), \quad \Pi(0,0)\in \BR^{n\times 2}, \quad S(0,0)=S(0,0)^*,  
 \\ &  \label{3.20}
 \I J  \in {\mathrm{GL}}(2, \BR), \quad J=J^*=J^{-1}, \quad  qJ=-J{ q}^*, \quad  QJ=-J  Q^*
 \end{align}
are valid,
where $q$ and $Q$ in \eqref{3.20} are given by \eqref{2.7}. Assume additionally that
\begin{align} &       \label{3.21}
d:=\det\Big(\big(I_2-\I J\Pi(0,0)^*S(0,0)^{-1}\cla^{-1}\Pi(0,0)\big)u(0,0)\Big)>0.
\end{align}
Then, the scalar function $\a$ and the matrix function $\wt u$ of the form
\begin{align} &       \label{3.22}
\wt u(\xi,\eta)=\a(\xi,\eta)d^{-1/2}\clu(\xi,\eta)u(\xi, \eta) ,
\end{align}
where $\clu(\xi,\eta)$ is given by \eqref{3.7}, satisfy \eqref{1.1}--\eqref{3.4+}.
\end{Cy}
\begin{proof}. According to \eqref{1.2} and \eqref{2.7}, we have
\begin{align} &       \label{3.23}
q(\xi,\eta), \, Q(\xi,\eta)\in \BR^{2\times 2}.
\end{align}
Relations \eqref{2.9}--\eqref{2.13}, \eqref{3.19}, \eqref{3.20} and \eqref{3.23}
show that
 \begin{align} &       \label{3.24}
A(\xi,\eta), \, S(\xi,\eta)\in \BR^{n\times n}, \quad \Pi(\xi,\eta)\in \BR^{n\times 2}.
 \end{align}
It follows from \eqref{1.2}, \eqref{3.7}, the first relation in \eqref{3.20}, and \eqref{3.24} that
\begin{align} &       \label{3.25}
\clu(\xi,\eta)\in \BR^{2\times 2}, \quad \wh u(\xi,\eta)\in \BR^{2\times 2}.
 \end{align}
 where $\wh u(\xi,\eta)=\clu(\xi,\eta)u(\xi,\eta)$.
Moreover, equalities \eqref{2.7} and the last two equalities in \eqref{3.20} imply that
$$\big( u(\xi,\eta)^*J u(\xi,\eta)\big)_{\xi}=0,  \quad \big( u(\xi,\eta)^*J u(\xi,\eta)\big)_{\eta}=0,$$
that is, 
\begin{align} &       \label{3.26}
 u(\xi,\eta)^*J u(\xi,\eta)= u(0,0)^*J u(0,0).
 \end{align}
 By virtue of \eqref{3.17} and \eqref{3.26} we have 
$\wh u(\xi,\eta)^*J\wh u(\xi,\eta) =u(0,0)^*J u(0,0)$.
Hence, taking taking into account that $\wh u(\xi,\eta)$ is continuous,
$\wh u(\xi,\eta)\in \BR^{2\times 2}$ and \eqref{3.21} holds we obtain
\begin{align} &       \label{3.27}
 \det \wh u(\xi,\eta)\equiv d>0.
 \end{align}
 According to Theorem \ref{TmSM}, $\a$ and $\wh u$ satisfy \eqref{1.1}.
 In view of \eqref{3.25} and \eqref{3.27}, $\wh u$ satisfies \eqref{1.2}
 and $\det \wh u(\xi,\eta)>0$. Now, compare \eqref{3.4} and \eqref{3.22} in order
 to see that $\a$ and $\wt u$ satisfy \eqref{1.1}--\eqref{3.4+}.
\end{proof}

\section{Ernst-type equations} \label{Ernst}
\setcounter{equation}{0}
{\bf 1.}  Non-isospectral (or modified) canonical system has the form 
\begin{align}& \label{E1}
w_{\xi}(\xi,z)=\I \la J H(\xi) w(\xi, z) \quad (\la=(z-\xi)^{-1}),
\\ &\label{E2}
H(\xi)=H(\xi)^*\in \BC^{m\times m}, \quad J=J^*=J^{-1}\in \BC^{m\times m},
\end{align}
where $H(\xi)\geq 0$. This system (or the corresponding multiplicative integrals) appeared, e.g.,  in the
works by M.S. Liv\v{s}ic} \cite{Liv}, by V.P. Potapov \cite{Pot}, by Yu.P. Ginzburg and by L.A. Sakhnovich (see the review \cite[pp. 37, 38]{SaLrev}).
It is closely connected (see \cite[pp. 34--39]{SaLrev}) with the Riemann-Hilbert problem for random matrices presented in \cite{D2}
and with Wiener-Masani problem in prediction theory (as discussed in \cite{SaLmn}).

GBDT for the canonical  system, that is, for system \eqref{E1}, \eqref{E2} ($H\geq 0$), where
the spectral parameter $\la$ does not depend on $z$ and $\xi$
was treated in \cite{SaA3}. GBDT for the non-isospectral system
\eqref{E1}, \eqref{E2} was studied in  \cite{ALS07}.

We note that B\"acklund transformation for Ernst equation was first introduced in \cite{Har}, and for the matrix form of Ernst equation
see, for instance, \cite{Alev, Ward}.
In particular, G.A. Alekseev \cite{Alev} considered Ernst equation as the compatibility condition
for the systems
\begin{align} &       \label{E3}
w_{\xi}=(z-\xi)^{-1}U(\xi, \eta)w, \quad w_{\eta}=(z-\eta)^{-1}V(\xi,\eta))w,
 \end{align}
where $U$ and $V$ have real-valued entries (for the hyperbolic case) and some special structure (see
\cite[(8),(11)]{Alev}). 

Somewhat modifying systems  \eqref{E3}, we consider the compatibility condition of the auxiliary linear systems
\begin{align} &       \label{E3'}
w_{\xi}=(z-\xi-\eta)^{-1}U(\xi, \eta)w, \quad w_{\eta}=(z-\xi- \eta)^{-1}V(\xi,\eta)w,
 \end{align}
and obtain an  Ernst-type  integrable nonlinear system (non-isospectral case):
\begin{align} &       \label{E4-}
U(\xi,\eta)-V(\xi,\eta)+[U(\xi,\eta),V(\xi,\eta)]=0, \quad U_{\eta}(\xi,\eta)=V_{\xi}(\xi,\eta).
\end{align}
Indeed, the compatibility condition \eqref{2.1} for systems \eqref{E3'} takes the form
\begin{align} & \nn
(z-\xi-\eta)^{-1}\big(U_{\eta}(\xi,\eta)-V_{\xi}(\xi,\eta)\big)
\\ &        \label{E4-0}
+(z-\xi-\eta)^{-2} \big(U(\xi,\eta)-V(\xi,\eta)+[U(\xi,\eta),V(\xi,\eta)]\big)=0,
\end{align}
which is equivalent to \eqref{E4-}. {\it We note that in the case of the systems
\eqref{E3'} it is convenient  to get rid of the spectral parameter $\la$ and use
the expression
$(z-\xi-\eta)^{-1}$ instead of it.} Thus, we deal with $w(\xi,\eta,z)$ where $z$
is the independent ``hidden" spectral parameter.

Further we set 
\begin{align} &       \label{E4}
U(\xi,\eta)= \I J H(\xi,\eta), \quad  V(\xi,\eta)=\I J \clh(\xi,\eta) \quad (\clh=\clh^*),
 \end{align}
and assume that \eqref{E2} holds. System \eqref{E4-} takes the form
\begin{equation}        \label{E4+}
JH(\xi,\eta)-J\clh(\xi,\eta)+\I [J H(\xi,\eta), J\clh(\xi,\eta)]=0, \quad H_{\eta}(\xi,\eta)=\clh_{\xi}(\xi,\eta).
\end{equation}

{\bf 2.}  In order to construct Darboux matrix corresponding to the system \eqref{E4+}, we  fix a triple
$\{\cla, S(0,0), \Pi(0,0)\}$ satisfying \eqref{2.8} and set
\begin{align} &       \label{E5}
A(\xi,\eta)= \big(\cla-(\xi+\eta)I_n\big)^{-1} \quad {\mathrm{i.e.}}, \quad A_{\xi}=A_{\eta}=A^2.
 \end{align}
 We introduce $\Pi(\xi,\eta)$ and $S(\xi,\eta)$ by the linear equations
\begin{align} &       \label{E6}
\Pi_{\xi}=- \I A \Pi J H, \quad \Pi_{\eta}=- \I A \Pi J \clh; \\
&       \label{E7-}
S_{\xi}= \Pi J H J^* \Pi^*- (AS+SA^*), \quad S_{\eta}= \Pi J \clh J^* \Pi^*- (AS+SA^*).
 \end{align}
 It is easily checked that by virtue of  \eqref{E4+},  \eqref{E6} and the last two equalities in \eqref{E5}
 we have $\Pi_{\xi \eta}=\Pi_{\eta \xi}$, that is, the  compatibility condition (for systems \eqref{E6}) is fulfilled.
 Moreover, the identity \eqref{2.15} is valid (see \cite[(2.4)]{ALS07}).

Now, we introduce a matrix function
\begin{align} &       \label{E7}
v(\xi,\eta,z):=w_0(\xi, \eta)^{-1}w_A\big(\xi,\eta, (z-\xi-\eta)^{-1}\big);
 \end{align}
where $w_A$ is given by \eqref{2.8+} and 
\begin{align} &       \label{E8}
\frac{\p}{\p \xi}w_0(\xi, \eta)= \wt G_0(\xi,\eta)  w_0(\xi,\eta), \quad \frac{\p}{\p \eta}w_0(\xi, \eta)= \wt F_0(\xi,\eta)  w_0(\xi,\eta),
\\ & \label{E9}
\wt G_0=-\I J\Pi^*S^{-1} \Pi-  [J \Pi^* S^{-1} \Pi, JH], \\
& \label{E10}
 \wt F_0=-\I J\Pi^*S^{-1} \Pi-  [J \Pi^* S^{-1} \Pi, J\clh],  \quad w_0(0,0)^*J w_0(0,0)=J.
 \end{align}
Using matrix functions  $A$, $\Pi$, $S$ and $w_0$ from above and taking into account our results for GBDT of the non-isospectral canonical system \cite{ALS07},
 we prove the following theorem.
 \begin{Tm}\label{TmpErnst} Let $H$ and $\clh$ satisfy \eqref{E4+}, let  the equalities \eqref{E2}  hold for $H$ and $J$ and assume that $\clh=\clh^*$.
Then, the matrix function $v(\xi,\eta,z)$ of the form \eqref{E7} is the corresponding Darboux matrix, that is, it satisfies the
systems:
 \begin{align} &       \label{E12}
v_{\xi}(\xi,\eta,z) =\I (z-\xi-\eta)^{-1}\big(J \wt H(\xi,\eta)v(\xi,\eta,z)-v(\xi,\eta,z)JH(\xi,\eta)\big),
\\ &       \label{E13}
v_{\eta}(\xi,\eta,z) =\I (z-\xi-\eta)^{-1}\big(J \wt \clh(\xi,\eta)v(\xi,\eta,z)-v(\xi,\eta,z)J \clh(\xi,\eta)\big),
 \end{align} 
where
\begin{align} &       \label{E14}
\wt H(\xi,\eta)=w_0(\xi, \eta)^*H(\xi,\eta)w_0(\xi,\eta), \quad \wt \clh(\xi,\eta)=w_0(\xi, \eta)^*\clh(\xi,\eta)w_0(\xi,\eta).
 \end{align} 
 \end{Tm}
 \begin{proof}. For each fixed $\xi$ or $\eta$ we substitute into \cite{ALS07} $z-\xi$ instead $z$ and $(\cla-(\xi+\eta) I_n)^{-1}$  instead of
 $A(\eta)$ or $z-\eta$ instead $z$ and $(\cla-(\xi+\eta) I_n)^{-1}$  instead of
 $A(\xi)$, respectively, and use \cite[(2.14), (2.15)]{ALS07} in order to derive
 \begin{align} &       \nn
\frac{d}{d \xi}w_A\big(\xi,\eta, (z-\xi-\eta)^{-1}\big)
\\ &  \label{E15}
=\I(z-\xi-\eta)^{-1}\big(JH(\xi,\eta)w_A\big(\xi,\eta, (z-\xi-\eta)^{-1}\big)
\\ & \quad \nn
-w_A\big(\xi,\eta, (z-\xi-\eta)^{-1}\big)JH(\xi,\eta)\big)+\wt G_0(\xi,\eta)w_A\big(\xi,\eta, (z-\xi-\eta)^{-1}\big)
 \end{align}
and
 \begin{align} &       \nn
\frac{d}{d \eta}w_A\big(\xi,\eta, (z-\xi-\eta)^{-1}\big)
\\ &  \label{E16}
=\I(z-\xi-\eta)^{-1}\big(J\clh(\xi,\eta)w_A\big(\xi,\eta, (z-\xi-\eta)^{-1}\big)
\\ & \quad \nn
-w_A\big(\xi,\eta, (z-\xi-\eta)^{-1}\big)J\clh(\xi,\eta)\big)+\wt F_0(\xi,\eta)w_A\big(\xi,\eta, (z-\xi-\eta)^{-1}\big),
 \end{align}
 where $\wt G_0$ and $\wt F_0$ are given by \eqref{E9} and \eqref{E10}.
 Here, we took into account that the definition \cite[(2.5)]{ALS07} of $w_A$ slightly differs from the
 definition in this paper. Using \eqref{E5}--\eqref{E7-} and \eqref{2.15}, one could also derive \eqref{E15} and \eqref{E16} directly.
 
 According to \eqref{E7}, \eqref{E8} and \eqref{E15}, \eqref{E16}, we have
 \eqref{E12} and \eqref{E13}, where
 \begin{align} &       \label{E17}
\wt H=Jw_0^{-1}JHw_0, \quad \wt \clh=Jw_0^{-1}\clh w_0 \quad(J=J^*=J^{-1}).
 \end{align} 
 In view of \eqref{E8}--\eqref{E10}, differentiating $w_0^*Jw_0$ we obtain
  \begin{align} &       \label{E18}
w_0(\xi,\eta)^*Jw_0(\xi, \eta)\equiv J.
 \end{align} 
 Finally, formulas \eqref{E17} and \eqref{E18} imply \eqref{E14}.
 \end{proof}
 \begin{Rk}\label{Rkw0} If $A$ is invertible, we may set
\begin{align} &       \label{E11}
w_0(\xi, \eta)=w_A(\xi,\eta,0)=I_m-\I J \Pi(\xi, \eta)^*S(\xi, \eta)^{-1}A(\xi, \eta)^{-1} \Pi(\xi, \eta).
 \end{align} 
 Indeed, in view of \cite[Remark 1]{ALS07} $w_A(\xi,\eta,0)$ satisfies \eqref{E8}.
Recall also that relation \eqref{2.15} yields  the identity \eqref{E19}, and so
the matrix $w_A(0,0,0)$ satisfies the last equality in \eqref{E10}.
 \end{Rk}
\begin{Rk}  According to \eqref{E17}, $J\wt H(\xi, \eta)$ is linear similar to $JH(\xi, \eta)$ and $J\wt \clh(\xi, \eta)$ is linear similar to $J \clh(\xi, \eta)$.
Moreover, in view of \eqref{E14} the inequality $H(\xi,\eta)\geq 0$ implies $\wt H(\xi,\eta)\geq 0$
and the inequality $\clh(\xi,\eta)\geq 0$ implies $\wt \clh(\xi,\eta)\geq 0$.
\end{Rk}

\section{Some examples} \label{Some}
\setcounter{equation}{0}
{\bf 1.} Explicit constructions are of special interest in our theory.
Recall that $\la$ of the form \eqref{3.8} satisfies \eqref{2.4} and \eqref{2.5}.
Compare \eqref{2.4} and \eqref{2.5} with \eqref{2.9} and \eqref{2.10}, respectively,
in order to see that the matrix function $A(\xi,\eta)$ in Theorem \ref{TmSM} may be
given explicitly. Namely, in view of Proposition \ref{PnQ} the following proposition is also valid.
\begin{Pn}\label{PnA} Let the matrix $\cla \in \BC^{n \times n}$ be given.
Assume that 
$$\la_k-2h(\eta)\not=0,  \quad \la_k+2f(\xi)\not=0$$ for the values $\la_k\in \s(\cla)$.
Then, the matrix function
\begin{align} &       \label{5.1}
A(\xi, \eta)=\Big(\clr\big(2h(\eta)\big)-\clr\big(-2f(\xi)\big)\Big)
\Big(\clr\big(2h(\eta)\big)+\clr\big(-2f(\xi)\big)\Big)^{-1},
\end{align} 
where $\clr$ is constructed in the proof of Proposition \ref{PnQ}, satisfies  \eqref{2.9} and \eqref{2.10}.
\end{Pn}
We note that $f, h$ and the square roots $\sqrt{\la_k-2h(\eta)}$ and $\sqrt{\la_k+2f(\xi)}$ in the construction
of $\clr$ should be chosen so that $\clr\big(2h(\eta)\big)+\clr\big(-2f(\xi)\big)$ is invertible and
$\sqrt{\la_k-2h(\eta)}$ and $\sqrt{\la_k+2f(\xi)}$ are
continuously differentiable. 

\begin{Rk}\label{Rku} Given $m=2p$ and $\a$ of the form \eqref{2.6}, \eqref{2.6'},
we may choose 
\begin{align} &       \label{3.5}
u(\xi,\eta)=\E^{\big(f(\xi)-h(\eta)\big)j},  \quad j=
\begin{bmatrix} I_p & 0 \\ 0 & -I_p
\end{bmatrix}, \end{align}
and
\begin{align} &       \label{3.5+}
J=\begin{bmatrix} 0 & I_p  \\  I_p & 0
\end{bmatrix} \quad {\mathrm{or}} \quad J=\begin{bmatrix} 0 &- \I I_p  \\ \I I_p & 0
\end{bmatrix}
\end{align}
in Theorem \ref{TmSM}. Indeed, in view of \eqref{2.7} and \eqref{3.5} we have
\begin{align} &       \label{3.6}
q(\xi,\eta)\equiv f^{\prime}(\xi)j, \quad Q(\xi,\eta)\equiv h^{\prime}(\eta)j \quad \Big(f^{\prime}=\frac{d}{d\xi}f\Big),
 \end{align} 
and the corresponding equalities in \eqref{2.14} hold. Clearly, \eqref{3.4!} holds as well. 
Finally, substituting \eqref{2.6} and \eqref{3.5} into the left-hand side of \eqref{1.1} we rewrite \eqref{1.1}
in the form $h^{\prime}(\eta)f^{\prime}(\xi)j-f^{\prime}(\xi)h^{\prime}(\eta)j=0$, and so \eqref{1.1} is valid.
In other words, $\a$ and $u$ given by  \eqref{2.6} and \eqref{3.5}, respectively, satisfy \eqref{1.1}.
\end{Rk}
Our next remark suggests the choice of $J$ and $u$ in Corollary \ref{Cygr}.
\begin{Rk} \label{Rkgr} Pauli matrix
$\s_2=J=\begin{bmatrix}0 &-\I \\ \I & 0 \end{bmatrix}$ gives a simple
example of the matrix $J$ such that the corresponding relations in \eqref{3.20} hold.

Setting $($in Remark \ref{Rku}$)$ $p=1$ and $J=\s_2$, we see that  relations \eqref{1.2} and \eqref{3.20}
are valid. According to Remark \ref{Rku}, the pair $\a,\, u$  satisfies \eqref{1.1}.
\end{Rk}

{\bf 2.}  Consider the case
 \begin{align} &       \label{5.2}
f(\xi)=-\xi, \quad h(\eta)=\eta,
 \end{align}  
 which was studied in \cite{MiYa}. According to \eqref{3.6} and \eqref{5.2} we have
 \begin{align} &       \label{5.2+}
q\equiv -j, \quad Q \equiv j.
 \end{align}  
 {\it New explicit solutions} appear when the parameter matrix $\cla$
 is non-diagonal.  In the next example, we deal with the simplest of such cases $(m=2p, \,\, n=2)$.
 \begin{Ee} Let $\cla$ be a $2\times 2$ Jordan block:
  \begin{align} &       \label{5.3}
\cla=\begin{bmatrix} c &1\\ 0 &c \end{bmatrix}.
 \end{align}  
 Then, in view of \eqref{B4} and \eqref{B5+} we have
 \begin{align} &       \label{5.4}
\clr(2\eta)=\begin{bmatrix} \om(\eta) & \frac{1}{2\om(\eta)} \\ 0 & \om(\eta) \end{bmatrix},
\quad  \clr(2\xi)=\begin{bmatrix} \nu(\xi) & \frac{1}{2\nu(\xi)} \\ 0 & \nu(\xi) \end{bmatrix}; \\
& \label{5.5}
 \om(\eta):=\sqrt{c-2\eta}, \quad
\nu(\xi):=\sqrt{c-2\xi}.
 \end{align} 
 After some simple calculations, using \eqref{5.1}, \eqref{5.2} and \eqref{5.4} we derive
 \begin{align}  &       \label{5.5+}
 A=\begin{bmatrix}a &b \\ 0 & a\end{bmatrix}, \quad a(\xi,\eta)=\frac{\om(\eta)-\nu(\xi)}{\om(\eta)+\nu(\xi)}, \quad 
 b(\xi,\eta)=-\frac{a(\xi,\eta)}{\nu(\xi)\om(\eta)}.
 \\  &       \label{5.6}
\big(A(\xi,\eta)-I_2\big)^{-1}=-\frac{\om(\eta)+\nu(\xi)}{2\nu(\xi)}\begin{bmatrix} 1 & \frac{\nu(\xi)- \om(\eta)}{2\om(\eta)\nu(\xi)^2} \\ 0 & 1 \end{bmatrix},
\\ & \label{5.7} 
\big(A(\xi,\eta)+I_2\big)^{-1}=\frac{\om(\eta)+\nu(\xi)}{2\om(\eta)}\begin{bmatrix} 1 & \frac{\om(\eta)-\nu(\xi)}{2\nu(\xi)\om(\eta)^2} \\ 0 & 1 \end{bmatrix}.
 \end{align} 
Partition $\Pi(\xi,\eta)$ into two $2 \times p$ blocks: $\Pi(\xi,\eta)=\begin{bmatrix} \Lam_1(\xi,\eta) & \Lam_2(\xi,\eta) \end{bmatrix}$,
recall that $\Pi(0,0)$ is assumed to be given $($it belongs to the triple, which determines GBDT$)$ and set 
\begin{align} &       \label{5.8}
\Lam_1(\xi,\eta)=\exp\left\{-\frac{1}{4}\begin{bmatrix} \big(\nu(\xi)+\om(\eta)\big)^2 &\frac{\nu(\xi)}{\om(\eta)}+\frac{\om(\eta)}{\nu(\xi)}\\ 0 &
 \big(\nu(\xi)+\om(\eta)\big)^2\end{bmatrix}\right\}\Lam_1(0,0),
\\
&       \label{5.9}
\Lam_2(\xi,\eta)=\exp\left\{\frac{1}{4}\begin{bmatrix} \big(\nu(\xi)+\om(\eta)\big)^2 &\frac{\nu(\xi)}{\om(\eta)}+\frac{\om(\eta)}{\nu(\xi)}\\ 0 &
 \big(\nu(\xi)+\om(\eta)\big)^2\end{bmatrix}\right\}\Lam_2(0,0),
 \end{align} 
where $\om$ and $\nu$ are introduced in \eqref{5.5}.  
Direct differentiation in   \eqref{5.8}, \eqref{5.9} and formulas \eqref{5.2+}, \eqref{5.6} and \eqref{5.7} show that we constructed
$\Pi(\xi,\eta)$ correctly and it satisfies \eqref{2.11}. Formulas \eqref{5.8} and \eqref{5.9} may be simplified
\begin{align} &       \label{5.10}
\Lam_1(\xi,\eta)=\exp\left\{- \big(\nu(\xi)+\om(\eta)\big)^2/4\right\}\left(I_2-\frac{1}{4}\begin{bmatrix} 0 &\frac{\nu(\xi)}{\om(\eta)}+\frac{\om(\eta)}{\nu(\xi)}\\ 0 & 0
\end{bmatrix}\right)
\Lam_1(0,0),
\\
&       \label{5.11}
\Lam_2(\xi,\eta)=\exp\left\{ \big(\nu(\xi)+\om(\eta)\big)^2/4\right\}\left(I_2+\frac{1}{4}\begin{bmatrix} 0 &\frac{\nu(\xi)}{\om(\eta)}+\frac{\om(\eta)}{\nu(\xi)}\\ 0 & 0
\end{bmatrix}\right)
\Lam_2(0,0),
 \end{align} 
 Finally, we note that under the condition
 $\om(\eta) \ov{\nu(\xi)}\not=\nu(\xi)\ov{\om(\eta)}$
 the  entries of the $2\times 2$  matrix function $S(\xi,\eta)=\{S_{ik}(\xi,\eta)\}_{i,k=1}^2$ are uniquely successively recovered from the identity \eqref{2.15}  $($and from formula \eqref{5.5+}$):$
 \begin{align} &      \nn
S_{22}=(a-\ov{a})^{-1}K_{22}, \quad S_{21}=(a-\ov{a})^{-1}(K_{21}+\ov{b}S_{22}), \\
& \nn
 S_{12}=(a-\ov{a})^{-1}(K_{12}-{b}S_{22}), \quad S_{11}=(a-\ov{a})^{-1}(K_{11}+\ov{b}S_{12}-bS_{21}),
 \end{align}  
where $K_{ik}$ are the entries of $K:=\I\Pi J\Pi^*$ and $\Pi(\xi,\eta)$ is explicitly constructed above.  Now, our main formulas
\eqref{3.7} and \eqref{3.7+} provide a corresponding family of explicit solutions $\wh u$ of the $\sigma$-model.
 \end{Ee}

{\bf Acknowledgments}  {This research    was supported by the
Austrian Science Fund (FWF) under Grant  No. P29177.}

\appendix
\section{Compatibility condition for  systems on $S$}\label{apA}
\setcounter{equation}{0}
Heuristically, the compatibility condition 
\begin{align} &       \label{A1}
S_{\xi \eta}=S_{\eta \xi}
 \end{align} 
for the systems  \eqref{2.12} and \eqref{2.13} may be deduced from the unique
solvability (under natural assumptions) of the identity \eqref{2.15} on $S$.

In order to  prove \eqref{A1} rigorously, we rewrite \eqref{2.9} and \eqref{2.10}
in the forms
\begin{align} &       \label{A2}
A_{\xi }=-\frac{\a_{\xi}}{\a}A(A+I_n)(A-I_n)^{-1},
\quad A_{\eta}=-\frac{\a_{\eta}}{\a}A(A-I_n)(A+I_n)^{-1}.
 \end{align} 
It is easy to see that the  identities 
\begin{align} &       \label{A3}
2A(A-I_n)^{-1}(A+I_n)^{-1}=(A-I_n)^{-1}+(A+I_n)^{-1},
\\  &       \label{A4}
2(A-I_n)^{-1}(A+I_n)^{-1}=(A-I_n)^{-1}-(A+I_n)^{-1}.
 \end{align} 
are valid. Now, we differentiate both sides of  \eqref{2.12} with respect to $\eta$ and both sides of \eqref{2.13}
with respect to $\xi$ using \eqref{2.13} and \eqref{2.12}, respectively, as well as the equalities \eqref{2.11} and \eqref{A2}.
We simplify the right-hand sides of the obtained relations using \eqref{2.14}, \eqref{A3} and \eqref{A4}.
Then, reducing similar terms we derive
\begin{align}       \nn
S_{\xi \eta}-S_{\eta \xi}=&-\frac{\I}{2\a}\Big(\a_{\xi}\big((A+I_n)^{-1}\Pi QJ\Pi^*(A^*+I_n)^{-1}
\\ & \nn
-(A-I_n)^{-1}\Pi QJ\Pi^*(A^*-I_n)^{-1}\big)
\\ & \nn
+\a_{\eta}\big((A+I_n)^{-1}\Pi qJ\Pi^*(A^*+I_n)^{-1}
\\ &  \label{A5}
-(A-I_n)^{-1}\Pi qJ\Pi^*(A^*-I_n)^{-1}\big)
\\ & \nn
+(A-I_n)^{-1}\Pi \big(2(\a q)_{\eta}+\a[Q,q]\big)J\Pi^*(A^*-I_n)^{-1}
\\ & \nn
-(A+I_n)^{-1}\Pi \big(2(\a Q)_{\xi}+\a[Q,q]\big)J\Pi^*(A^*+I_n)^{-1}\Big).
\end{align} 
Next, we multiply \eqref{A5} by $2\I \a (A-I_n)(A+I_n)$ from the left and by \\
$(A^*+I_n)(A^*-I_n)$ from the right. Taking into account that $(\a q)_{\eta}=(\a Q)_{\xi}$, we obtain
\begin{align}  &     \nn
2\I \a (A-I_n)(A+I_n)\big(S_{\xi \eta}-S_{\eta \xi}\big)(A^*+I_n)(A^*-I_n)
\\ &  \label{A6}
=-2\a_{\xi}\big(A\Pi QJ\Pi^*+\Pi QJ\Pi^*A^*\big)- 2\a_{\eta}\big(A\Pi qJ\Pi^*+\Pi qJ\Pi^*A^*\big)
\\ & \nn
+2A\Pi\big((\a q)_{\eta}+(\a Q)_{\xi}+\a[Q,q]\big)J\Pi^*
+2\Pi \big((\a q)_{\eta}+(\a Q)_{\xi}+\a[Q,q]\big)J\Pi^*A^*.
\end{align} 
Finally, the first equation in \eqref{C} implies that the right-hand side of \eqref{A6}
equals zero. Thus, \eqref{A1} follows.

\section{Matrix square roots}\label{apB}
\setcounter{equation}{0}
Modifying the proof of  \cite[Proposition 3.3]{ALS20}, we obtain the following proposition.
\begin{Pn}\label{PnQ}  Let $\cla\in \BC^{n\times n}$ admit representation
$E\clj E^{-1}$, where $\clj$ is the Jordan normal form of $\cla$. Then, there is a matrix function
$\clr(\mu)$ $(\mu\in \BR, \,\, \mu \not\in \s(\cla))$ such that 
\begin{align}& \label{B1}
\clr(\mu)^2=\cla -\mu I_n, \quad  \clr(\mu)=E\cld(\mu)E^{-1},
\\ & \label{B2}
\clr(\mu_1)\clr(\mu_2)=\clr(\mu_2)\clr(\mu_1) \quad (\mu_1,\mu_2\in \BR),
\end{align}
where $\cld$ is a block diagonal matrix with the blocks of the same orders as the corresponding Jordan blocks
of $\clj$. Moreover, the blocks of $\cld$ are upper triangular Toeplitz matrices $($or scalars if the corresponding
blocks of $\clj$ are scalars$)$.
\end{Pn}
\begin{proof}. 
Clearly, the statement of proposition holds for $n=1$. Consider the case, where $\cla$ is an $n\times n$ Jordan block $(n\geq 2)$:
\begin{align}& \label{B2+}
\cla=\begin{bmatrix} \la & 1 & & 
\\  & \la & \ddots & \\
& & \ddots & 1 \\
& & & \la
\end{bmatrix}.
\end{align}

For this $\cla$ and $\mu\in \BR$, we construct  upper triangular Toeplitz matrices $\clr(\mu)$ satisfying the first equality in \eqref{B1}.
First, we introduce the shift matrices
\begin{align}& \label{B3}
 \cls_{i}:=\{\delta_{k-l+i}\}_{k,l=1}^n, \quad \cls_i \cls_j =\cls_{i+j}.
\end{align}
where $\delta_s$ is Kronecker delta, and $\cls_{i}=0$ for $i \geq n$.
Let us write down the representations
\begin{align}& \label{B4}
\cla-\mu I_n=(\la-\mu) I_n+\cls_1, \quad \clr(\mu)=c_0 I_n +c_1 \cls_1 + \ldots + c_{n-1}\cls_{n-1}.
\end{align}
According to \eqref{B3} and \eqref{B4}, we have
\begin{align} \nn
\clr(\mu)^2=&c_0^2I_n+2c_0c_1\cls_1+(2c_0c_2+c_1^2)\cls_2
\\ & \label{B5}
+ \sum_{i=3}^{n-1}\big(2c_0 c_i+c_1c_{i-1}+\ldots +c_{i-1}c_1\big)\cls_i.
\end{align}
Now, we set 
\begin{align}& \label{B5+}
c_0=\sqrt{\la-\mu}\not=0,  \quad c_1=1/(2c_0),
\end{align}
and choose successively the values
$c_2, \, \ldots$ so that the coefficients before the shift matrices $\cls_i$ $(i\geq 2)$ on the right-hand side
of \eqref{B5} turn to zero. (In this way, the values $c_i$ $(i\geq 2)$ are uniquely defined, and the upper triangular Toeplitz matrix $\clr(\mu)$ 
given by \eqref{B4} satisfies the first equality in \eqref{B1}.)

When $\cla$ is a Jordan matrix $\clj$,  we construct block diagonal  matrix $\cld(\mu)$, each block of which is generated
by the corresponding Jordan block of $\clj$ in a way described above. It is easy to see that
\eqref{B1} holds for $\cla=\clj $ and $\clr(\mu)=\cld(\mu)$. Finally, if $\cla=E\clj E^{-1}$, we set $\clr(\mu)=E \cld(\mu) E^{-1}$
(as in the second equality in \eqref{B1}),
and the first equality in \eqref{B1} for $\cla$ and $\clr$ follows from  the first equality in \eqref{B1} for $\clj$ and $\cld$.

Since the blocks of $\clj$ and $\cld(\mu)$ are upper triangular  Toeplitz matrices, $\clj$ and $\cld(\mu)$ commute
(see, e.g., \cite{Com} on the properties of triangular Toeplitz matrices). Hence, \eqref{B2} is valid.
\end{proof}
Although $\mu \in \BR$ is required is Section \ref{Some}, it is easy to see that the construction above
works also for  $\mu \in \BC$.

\begin{flushright}
A.L. Sakhnovich,\\
Faculty of Mathematics,
University
of
Vienna, \\
Oskar-Morgenstern-Platz 1, A-1090 Vienna,
Austria, \\
e-mail: {\tt oleksandr.sakhnovych@univie.ac.at}

\end{flushright}

\end{document}